\documentclass[12pt]{article}
\usepackage{latexsym,amsfonts,amssymb}
\usepackage[dvips]{color}
\usepackage{colordvi,multicol}
\usepackage{amsmath}
\usepackage{graphicx}
\usepackage{pdflscape}
\usepackage{rotating}

\usepackage{latexsym,amsfonts,amssymb}
\usepackage[dvips]{color}
\usepackage{colordvi,multicol}
\usepackage{amsmath}
\usepackage{graphicx}
\usepackage{pdflscape}
\usepackage{rotating}
\usepackage{breqn}

\newtheorem{thm}{Theorem}[section]

\newtheorem{rem}[thm]{Remark}

\newtheorem{que}[thm]{Question}

\date{}
\begin{document}

\title{\bf The extreme values of two probability functions for the Gamma distribution}
\author{Ping Sun, Ze-Chun Hu  and Wei Sun\thanks{Corresponding author.}\\ \\
  {\small Business School, Chengdu University, Chengdu 610106, China}\\ \\
 {\small College of Mathematics, Sichuan University, Chengdu 610065, China}\\ \\
{\small Department of Mathematics and Statistics, Concordia University, Montreal H3G 1M8,  Canada}\\ \\
{\small  sunping@cdu.edu.cn\ \ \ \ zchu@scu.edu.cn\ \ \ \ wei.sun@concordia.ca}}

\maketitle

\begin{abstract}
\noindent
Motivated by Chv\'{a}tal's conjecture and Tomaszewaki's conjecture, we investigate
the extreme value problem of two probability functions for the Gamma distribution.
Let $\alpha,\beta$ be arbitrary positive real numbers and $X_{\alpha,\beta}$ be a Gamma random variable with shape parameter $\alpha$ and scale parameter $\beta$.
We study the extreme values of functions  $P\{X_{\alpha,\beta}\le E[X_{\alpha,\beta}]\}$ and $P\{|X_{\alpha,\beta}-E[X_{\alpha,\beta}]|\le \sqrt{{\rm Var}(X_{\alpha,\beta})}\}$. Among other things, we show that $
\inf_{\alpha,\beta}P\{X_{\alpha,\beta}\le E[X_{\alpha,\beta}]\}=\frac{1}{2}$ and $\inf_{\alpha,\beta}P\{|X_{\alpha,\beta}-E[X_{\alpha,\beta}]|\le \sqrt{{\rm Var}(X_{\alpha,\beta})}\}=P\{|Z|\le 1\}\approx 0.6826$, where $Z$ is a standard normal random variable.
\end{abstract}



\section{Introduction}

Special probability distributions play a fundamental role in probability theory, statistics, optimization and different research fields of science including physics, chemistry, bioscience, economy and management science, etc. Many famous distributions have a long history and date back to the era of   Bernoulli and De
Moivre.  However, our understanding of them is far from complete. This paper is motivated by Chv\'{a}tal's conjecture for the binomial distribution and Tomaszewaki's conjecture for the  Rademacher
sequence,
both of which were completely solved very recently.

Let $B(n,p)$ denote a binomial random variable with parameters $n$ and $p$. Janson \cite{Ja21} introduced the following conjecture suggested by Va\v{s}k Chv\'{a}tal.

\noindent {\bf Conjecture 1} (Chv\'{a}tal)\ \  For any fixed $n\geq 2$, as $m$ ranges over $\{0,\ldots,n\}$, the probability $q_m:=P\{B(n,\frac{m}{n})\leq m\}$ is the smallest when
 $m$ is the integer closest to $\frac{2n}{3}$.

\noindent   Chv\'{a}tal's conjecture has interesting applications in machine learning. Janson \cite{Ja21} proved Chv\'{a}tal's conjecture for sufficiently large $n$. Barabesi et al. \cite{BPR23} and Sun \cite{Su21} showed that Chv\'{a}tal's conjecture is true for any $n\geq 2$.

The second motivation of this paper is from the problem attributed to Boguslav Tomaszewski.

\noindent {\bf Conjecture 2} (Tomaszewski)\ \  Let $X=\sum_{i=1}^na_ix_i$, where $\sum_{i=1}^n a_i^2=1$ and  $\{x_i\}$ is a sequence of independent $\{-1,1\}$-valued symmetric random variables. Then, $P\{|X|\leq 1\}\geq 1/2$.

\noindent Tomaszewski's conjecture has many applications in probability theory, geometric analysis and computer science. Recently Keller and Klein \cite{KK22} completely solved Tomaszewski's conjecture. We refer the reader to Keller and  Klein \cite{KK22} for the details, and Dvorak and Klein \cite{DK22}  and Hu et al. \cite{Hu} for some related problems.

In this paper, we will focus on the Gamma distribution. It is well-known that the Gamma distribution, including  the exponential and the $\chi^2$-distributions as two important special cases, is one of the most basic probability distributions. It has a nearly 200-year history dating back to Laplace.
The Gamma distribution has been widely used in practice. It is frequently applied to describe the time between independent events that occur at a constant average rate. According to Wikipedia \cite{Wiki}, there are many significant applications. For example, it has been used to model the size of insurance claims and rainfalls, the multi-path fading of signal power, the age distribution of cancer incidence, the copy number of a constitutively expressed protein, etc.

Motivated by Chv\'{a}tal's conjecture,  Xu et al. \cite{XLH22} initiated the study of the minimum value problem for special probability distributions. Let $\{Y_{\lambda}\}$ be a family of random variables with the same distribution $F$ but  different parameters $\lambda$. Define
$$
r(\lambda):=P\{Y_{\lambda}\le E[Y_{\lambda}]\}.
$$
 Xu et al. \cite{XLH22} discussed the minimum value of the function $r$ and gave a complete answer if $F$ is the Poisson distribution or the geometric distribution. Further, Li and Hu \cite{LH23} considered the  minimum value problem for the $\chi^2$, Weibull and Pareto distributions.

In the first part of this paper, we consider the following more general minimum value problem for the Gamma distribution. Let $\alpha,\beta,\kappa>0$ and $X_{\alpha,\beta}$ be a Gamma random variable with shape parameter $\alpha$ and scale parameter $\beta$.
Define
\begin{equation}\label{Mar29a}
g_{\kappa}(\alpha,\beta):=P\{X_{\alpha,\beta}\le \kappa E[X_{\alpha,\beta}]\}.
\end{equation}
For fixed $\kappa$, what is the minimum value of the function $g_{\kappa}$? In Section 2, we will give a complete answer to this question. Interestingly, we discovered an unnoticed phase transition phenomenon (cf. Figures 1-4 and Remark \ref{rem26}) and obtained the following
result.
\begin{thm}\label{pro22}
Let $\alpha,\beta$ be arbitrary positive real numbers and $X_{\alpha,\beta}$  be a Gamma random variable with shape parameter $\alpha$ and scale parameter $\beta$. Then,
\begin{equation}\label{inf2}
P\{X_{\alpha,\beta}\le E[X_{\alpha,\beta}]\}>\frac{1}{2},
\end{equation}
and
$$
\inf_{\alpha,\beta}P\{X_{\alpha,\beta}\le E[X_{\alpha,\beta}]\}=\frac{1}{2}.
$$
\end{thm}

In the second part of this paper, we prove a more interesting and deeper 
result.
\begin{thm}\label{thm11}
Let $\alpha,\beta$ be arbitrary positive real numbers,  $X_{\alpha,\beta}$  be a Gamma random variable with shape parameter $\alpha$ and scale parameter $\beta$, and $Z$ be a standard normal random variable. Then,
\begin{equation}\label{inf}
P\{|X_{\alpha,\beta}-E[X_{\alpha,\beta}]|\le \sqrt{{\rm Var}(X_{\alpha,\beta})}\}>P\{|Z|\le 1\}\approx 0.6826,
\end{equation}
and
$$
\inf_{\alpha,\beta}P\{|X_{\alpha,\beta}-E[X_{\alpha,\beta}]|\le \sqrt{{\rm Var}(X_{\alpha,\beta})}\}=P\{|Z|\le 1\}.
$$
\end{thm}

The proof of inequality (\ref{inf}) will be given in Sections 3 and 4. It is worth pointing out that the software {\bf Mathematica} plays a crucial role in our work, which provides us with deep insight into how to handle delicate inequalities. However, all proofs contained in
this paper remain rigorous and easily verifiable.

Note that if $0<\kappa\not=1$, then the following more general inequality might not hold for all Gamma random variables:
$$
P\{|X_{\alpha,\beta}-E[X_{\alpha,\beta}]|\le \kappa\sqrt{{\rm Var}(X_{\alpha,\beta})}\}>P\{|Z|\le \kappa\}.
$$
For example, we have
\begin{eqnarray*}
&& 0.3834005=P\{|X_{1,1}-E[X_{1,1}]|\le 0.5\sqrt{{\rm Var}(X_{1,1})}\}\\
&>&0.3829249=P\{|Z|\le 0.5\}\\
&>&0.3819693=P\{|X_{2,1}-E[X_{2,1}]|\le 0.5\sqrt{{\rm Var}(X_{2,1})}\},
\end{eqnarray*}
and
\begin{eqnarray*}
&&0.9502129=P\{|X_{1,1}-E[X_{1,1}]|\le 2\sqrt{{\rm Var}(X_{1,1})}\}\\
&<&0.9544997=P\{|Z|\le 2\}\\
&<&0.9585112=P\{|X_{10,1}-E[X_{10,1}]|\le 2\sqrt{{\rm Var}(X_{10,1})}\}.
\end{eqnarray*}

We know that the Gamma distribution is infinitely divisible and, for any L\'evy process $\{L_t,t\ge0\}$, the distribution of $L_t$   is infinitely divisible. Motivated by Theorem \ref{thm11}, it is natural to ask if any infinitely divisible  random variable $L$  satisfies the following inequality:
\begin{equation}\label{inf30}
P\{|L-E[L]|\le \sqrt{{\rm Var}(L)}\}\ge P\{|Z|\le 1\}.
\end{equation}
We have done  numerical calculations, but not yet found a familiar  infinitely divisible random variable which does not satisfy (\ref{inf30}). We pose the following question:
\begin{que}
Does inequality (\ref{inf30}) hold for any  infinitely divisible random variable $L$? If yes, give a proof. If no, can we give good sufficient conditions which ensure the validity of (\ref{inf30})?
\end{que}

\section{The minimum value problem for the function $g_{\kappa}$}\setcounter{equation}{0}

Let $\alpha,\beta>0$ and $X_{\alpha,\beta}$ be a Gamma random variable with probability density function:
$$
f_{\alpha,\beta}(x)=\frac{x^{\alpha-1}e^{-x/\beta}}{\Gamma(\alpha)\beta^{\alpha}},\ \ \ \ x>0.
$$
For $\kappa>0$, we consider the extreme values of the probability function  $g_{\kappa}$ defined by (\ref{Mar29a}).

We have
\begin{eqnarray*}
g_{\kappa}(\alpha,\beta)=\int_0^{\kappa\alpha\beta}\frac{x^{\alpha-1}e^{-x/\beta}}{\Gamma(\alpha)\beta^{\alpha}}dx=
\int_0^{\kappa\alpha}\frac{y^{\alpha-1}e^{-y}}{\Gamma(\alpha)}dy=g_{\kappa}(\alpha,1).
\end{eqnarray*}
Then, we may assume without loss of generality that $\beta=1$ and focus on the extreme values of the following function
$$
h_{\kappa}(\alpha):=g_{\kappa}(\alpha,1)=\int_0^{\kappa\alpha}\frac{y^{\alpha-1}e^{-y}}{\Gamma(\alpha)}dy,\ \ \ \ \alpha>0.
$$

By Euler's reflection formula
$$
\Gamma(1-\alpha)\Gamma(\alpha)=\frac{\pi}{\sin(\pi\alpha)},\ \ \ \ \alpha\in(0,1),
$$
we get
\begin{eqnarray}\label{25a}
\liminf_{\alpha\downarrow0}h_{\kappa}(\alpha)
&=&\liminf_{\alpha\downarrow0}\int_0^{\kappa\alpha}\alpha y^{\alpha-1}e^{-y}dy\nonumber\\
&\ge&\liminf_{\alpha\downarrow0}\left[\alpha e^{-\kappa\alpha}\int_0^{\kappa\alpha} y^{\alpha-1}dy\right]\nonumber\\
&=&\liminf_{\alpha\downarrow0}(\kappa\alpha)^{\alpha} \nonumber\\
&=&\liminf_{\alpha\downarrow0}e^{\alpha\ln(\kappa\alpha)} \nonumber\\
&=&1.
\end{eqnarray}
It follows that for any $\kappa>0$,
$$
\sup_{\alpha>0}h_{\kappa}(\alpha)=1.
$$

In the sequel, we only consider the infimum value of  the function $h_{\kappa}$.

\subsection{Case $\kappa\le 1$}

In this subsection, we assume that $\kappa\le 1$. We will show that $h_{\kappa}(\alpha+1)<h_{\kappa}(\alpha)$ for any $\alpha>0$. In fact, we have that
\begin{eqnarray}\label{22a}
h_{\kappa}(\alpha+1)<h_{\kappa}(\alpha)
&\Leftrightarrow&\Gamma(\alpha+1)[h_{\kappa}(\alpha+1)-h_{\kappa}(\alpha)]<0\nonumber\\
&\Leftrightarrow&\int_0^{\kappa(\alpha+1)}y^{\alpha}e^{-y}dy-\alpha\int_0^{\kappa\alpha}y^{\alpha-1}e^{-y}dy<0\nonumber\\
&\Leftrightarrow&-[\kappa(\alpha+1)]^{\alpha}e^{-\kappa(\alpha+1)}+\alpha\int_{\kappa\alpha}^{\kappa(\alpha+1)}y^{\alpha-1}e^{-y}dy<0\nonumber\\
&\Leftrightarrow&\alpha\int_{\kappa\alpha}^{\kappa(\alpha+1)}y^{\alpha-1}e^{\kappa(\alpha+1)-y}dy<[\kappa(\alpha+1)]^{\alpha}\nonumber\\
&\Leftrightarrow&\alpha\int_{0}^{\kappa}[\kappa(\alpha+1)-w]^{\alpha-1}e^{w}dw<[\kappa(\alpha+1)]^{\alpha}\nonumber\\
&\Leftrightarrow&-\int_{0}^{\kappa}e^{w}d[\kappa(\alpha+1)-w]^{\alpha}<[\kappa(\alpha+1)]^{\alpha}\nonumber\\
&\Leftrightarrow&\int_{0}^{\kappa}[\kappa(\alpha+1)-w]^{\alpha}e^{w}dw<(\kappa\alpha)^{\alpha}e^{\kappa}\nonumber\\
&\Leftrightarrow&\int_{0}^{\kappa}\left(1+\frac{\kappa-w}{\kappa\alpha}\right)^{\alpha}e^{w-\kappa}dw<1\nonumber\\
&\Leftrightarrow&\int_{0}^{\kappa}\left(1+\frac{z}{\kappa\alpha}\right)^{\alpha}e^{-z}dz<1\nonumber\\
&\Leftrightarrow&\int_{0}^{1}\kappa\left(1+\frac{w}{\alpha}\right)^{\alpha}e^{-\kappa w}dw<1.
\end{eqnarray}
Note that $(1+\frac{1}{x})^x$ is strictly increasing with respect to $x>0$ and the limit equals $e$. Then, we have
$$
\left(1+\frac{w}{\alpha}\right)^{\alpha}=\left\{\left(1+\frac{w}{\alpha}\right)^{\frac{\alpha}{w}}\right\}^{w}<e^{w}.
$$
Thus,
\begin{eqnarray}\label{22b}
\int_{0}^{1}\kappa\left(1+\frac{w}{\alpha}\right)^{\alpha}e^{-\kappa w}dw
<\int_{0}^{1}e^{w}\cdot \kappa e^{-\kappa w}dw\le \int_{0}^{1}e^{w}\cdot e^{- w}dw=1.
\end{eqnarray}
Hence, by (\ref{22a}) and (\ref{22b}), we obtain that $\{h_{\kappa}(\alpha+n)\}_{n\in\mathbb{N}}$ is a strictly decreasing sequence.

Let $Y_0,Y_1,Y_2,\dots$ be independent random variables such that $Y_0\sim$ Gamma$(\alpha,1)$ and $Y_i\sim$ Gamma$(1,1)$, $i\ge 1$. Then, for $\kappa\in(0,1)$, by the strong law of large numbers, we get
\begin{eqnarray*}
\lim_{n\rightarrow\infty}h_{\kappa}(\alpha+n)
&=&\lim_{n\rightarrow\infty}P\{Y_0+Y_1+\cdots +Y_n\le \kappa E[Y_0+Y_1+\cdots +Y_n]\}\\
&=&\lim_{n\rightarrow\infty}P\left\{\frac{Y_0}{\alpha+n}+\frac{Y_1+\cdots +Y_n}{\alpha+n}\le \kappa \right\}\\
&=&0;
\end{eqnarray*}
and for $\kappa=1$, by the central limit theorem, we get
\begin{eqnarray*}
\lim_{n\rightarrow\infty}h_{\kappa}(\alpha+n)
&=&\lim_{n\rightarrow\infty}P\{Y_0+Y_1+\cdots +Y_n\le E[Y_0+Y_1+\cdots +Y_n]\}\\
&=&\lim_{n\rightarrow\infty}P\{(Y_0-\alpha)+(Y_1-1)+\cdots+(Y_n-1)\le 0\}\\
&=&\lim_{n\rightarrow\infty}P\left\{\frac{Y_0-\alpha}{\sqrt{n}}+\frac{(Y_1-1)+\cdots+(Y_n-1)}{\sqrt{n}}\le 0\right\}\\
&=&\frac{1}{2}.
\end{eqnarray*}
Therefore, 
\begin{eqnarray*}
\inf_{\alpha>0}h_{\kappa}(\alpha)=\left\{\begin{array}{ll}0,\ \ \ \  & \kappa\in(0,1), \\ \frac{1}{2},\ \ \ \  &
\kappa=1,\end{array}\right.\end{eqnarray*}
which implies  that inequality (\ref{inf2}) holds. The proof of Theorem \ref{pro22} is complete.

\begin{rem}  It is interesting to compare (\ref{inf2}) with the well-known fact that the Gamma distribution is skewed to the right and its skewness is given by
$$
E\left[\left(\frac{X_{\alpha,\beta}-E[X_{\alpha,\beta}]}{\sqrt{{\rm Var}(X_{\alpha,\beta})}}\right)^3\right]=\frac{2}{\sqrt{\alpha}}.
$$
\end{rem}

\subsection{Case $\kappa> 1$}

In this subsection, we assume that $\kappa> 1$. This case is especially interesting. By (\ref{inf2}), for any $\alpha>0$, we have
\begin{eqnarray}\label{25c}
h_{\kappa}(\alpha)=P\{X_{\alpha,1}\le \kappa E[X_{\alpha,1}]\}>P\{X_{\alpha,1}\le E[X_{\alpha,1}]\}>\frac{1}{2}.
\end{eqnarray}
Let $Y_1,Y_2,\dots$ be
independent Gamma$(1,1)$ random variables. By the central limit theorem, we get
\begin{eqnarray}\label{25d}
\liminf_{\alpha\rightarrow\infty}h_{\kappa}(\alpha)\nonumber
&\ge&\liminf_{\alpha\rightarrow\infty}P\{Y_1+\cdots +Y_{[\alpha]+1}\le \kappa[\alpha]\}\nonumber\\
&=&\lim_{n\rightarrow\infty}P\left\{\frac{(Y_1-1)+\cdots+(Y_{n+1}-1)}{\sqrt{n+1}}\le \frac{(\kappa-1)n-1}{\sqrt{n+1}}\right\}\nonumber\\
&=&1.
\end{eqnarray}
Hence, by (\ref{25a}), (\ref{25c}) and (\ref{25d}), we obtain that
$$
\min_{\alpha>0}h_{\kappa}(\alpha)>\frac{1}{2},\ \ \ \ \forall\kappa>1.
$$

Further, by virtue of {\bf Mathematica}, we obtain the following numerical results.
$$\min_{\alpha>0}h_{\kappa}(\alpha)=\left\{\begin{array}{ll}h_{\kappa}(33.4871)=0.545885,\ \ & \kappa=1.01, \\h_{\kappa}(3.47146)=0.64021,\ \ & \kappa=1.1, \\
h_{\kappa}(1.78959)=0.691283,\ \ & \kappa=1.2, \\
h_{\kappa}(0.757559)=0.774739,\ \ & \kappa=1.5, \\
h_{\kappa}(0.396184)=0.841243,\ \ & \kappa=2, \\ h_{\kappa}(0.205464)=0.899108, &
\kappa=3,\\
h_{\kappa}(0.13917)=0.925864,\ \ &
\kappa=4.\end{array}\right.$$
Below are graphs of the function $h_{\kappa}(\alpha)$ for different values of $\kappa$.

\begin{figure}[h]
\begin{center}
\scalebox{0.8}{\includegraphics{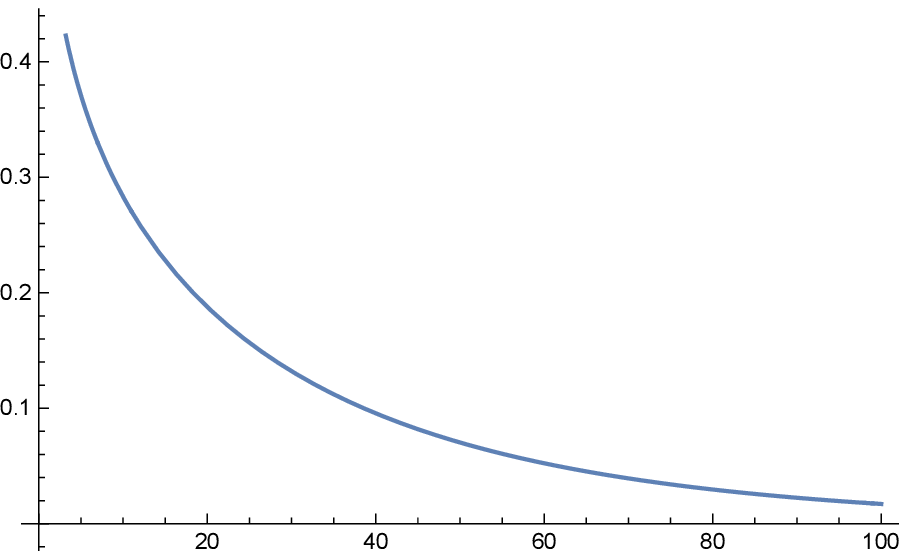}}
\end{center}
\end{figure}
\begin{center}
{\small Figure 1: Function $h_{\kappa}(\alpha)$: $\kappa=0.8$.}
\end{center}

\begin{figure}[h]
\begin{center}
\scalebox{0.8}{\includegraphics{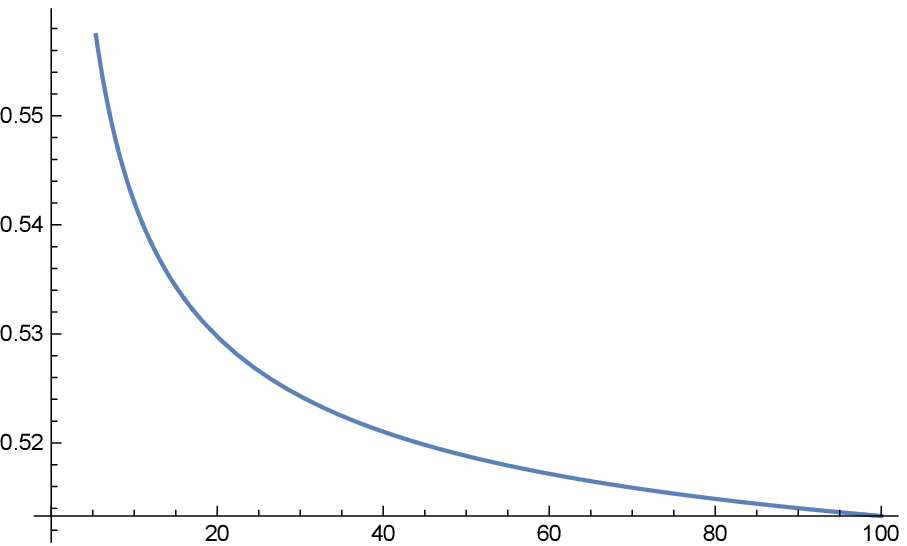}}
\end{center}
\end{figure}
\begin{center}
{\small Figure 2: Function $h_{\kappa}(\alpha)$: $\kappa=1$.}
\end{center}

\newpage\begin{figure}[h]
\begin{center}
\scalebox{0.8}{\includegraphics{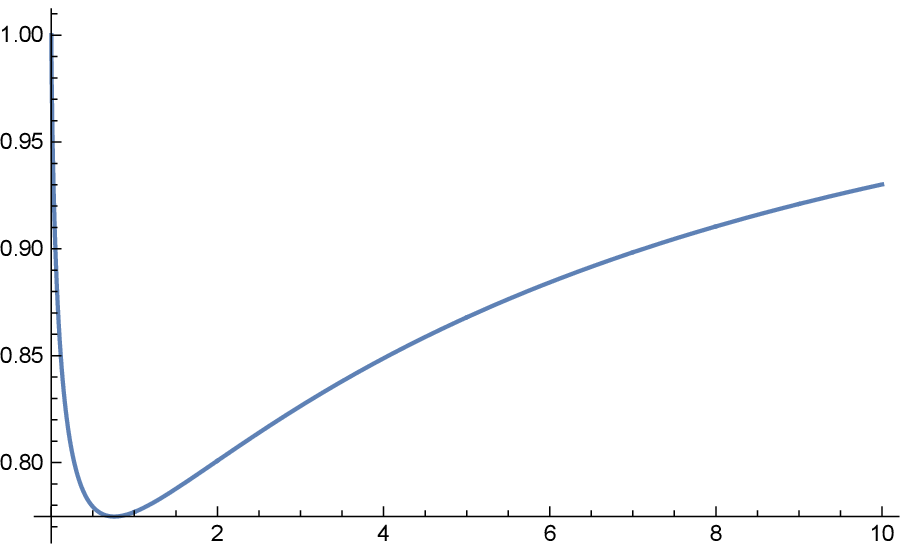}}
\end{center}
\end{figure}
\begin{center}
{\small Figure 3: Function $h_{\kappa}(\alpha)$: $\kappa=1.5$.}
\end{center}

\begin{figure}[h]
\begin{center}
\scalebox{0.8}{\includegraphics{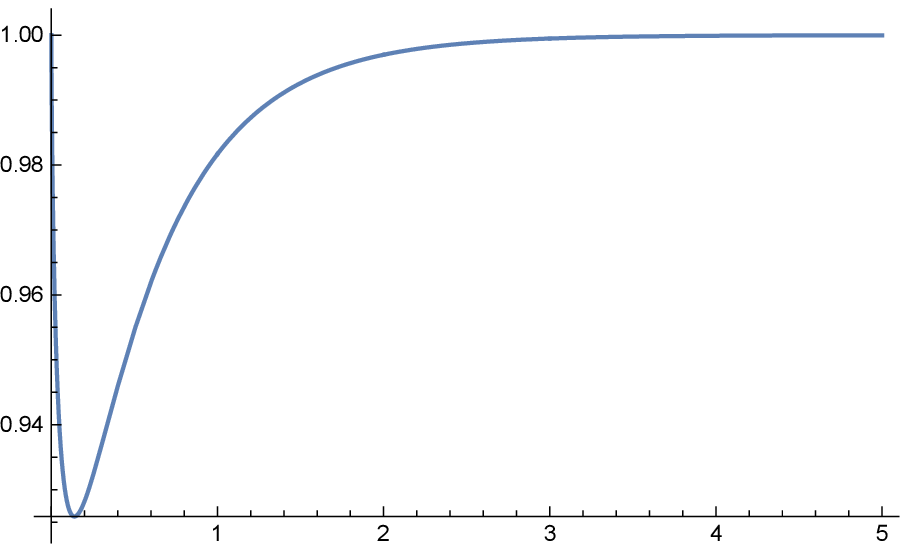}}
\end{center}
\end{figure}
\begin{center}
{\small Figure 4: Function $h_{\kappa}(\alpha)$: $\kappa=4$.}
\end{center}

\begin{rem}\label{rem26}
The above analysis
shows that there is an interesting phase transition phenomenon in the minimum value problem for the Gamma distribution. The critical point is $\kappa=1$ and the behaviors for the three cases $\kappa<1$, $\kappa=1$ and $\kappa>1$ are totally different.
\end{rem}

\section{Variation comparison between Gamma distribution and normal distribution}\setcounter{equation}{0}

Let $\alpha,\beta>0$. Define
$$
s(\alpha,\beta):=P\{|X_{\alpha,\beta}-E[X_{\alpha,\beta}]|\le \sqrt{{\rm Var}(X_{\alpha,\beta})}\}.
$$
We have
\begin{eqnarray*}
s(\alpha,\beta)=P\{|X_{\alpha,\beta}-\alpha\beta|\le \sqrt{\alpha}\beta\}=P\left\{\left|\frac{X_{\alpha,\beta}}{\beta}-\alpha\right|\le \sqrt{\alpha}\right\}=s(\alpha,1).
\end{eqnarray*}
Then, we may assume without loss of generality that $\beta=1$ and focus on  the following function
$$
t(\alpha):=s(\alpha,1)=\int_{\max\{0,\alpha-\alpha^{\frac{1}{2}}\}}^{\alpha+\alpha^{\frac{1}{2}}}\frac{y^{\alpha-1}e^{-y}}{\Gamma(\alpha)}dy,\ \ \ \ \alpha>0.
$$

 In this section, we will prove the following result.
\begin{thm}\label{Mar26v}
For  any $\alpha>0$,
\begin{equation}\label{Mar18a}
t(\alpha+1)<t(\alpha).
\end{equation}
\end{thm}

Note that by the central limit theorem,
$\lim_{\alpha\rightarrow\infty}t(\alpha)=P\{|Z|\leq 1\}$, where $Z$ is a standard normal random variable.
Once (\ref{Mar18a}) is proved, we conclude that
$$
s(\alpha,\beta)=s(\alpha,1)=t(\alpha)>P\{|Z|\leq 1\}\approx 0.6826,\ \ \ \ \forall \alpha,\beta>0.
$$
Then,  inequality (\ref{inf}) holds and hence the proof of Theorem \ref{thm11} is complete.
\vskip 0.3cm
\noindent {\bf Proof of Theorem \ref{Mar26v}.}\ \ First, we consider the case $0<\alpha\le1$. We have that
\begin{eqnarray}\label{Mar26m}
&&t(\alpha+1)<t(\alpha)\nonumber\\
&\Leftrightarrow&\Gamma(\alpha+1)[t(\alpha+1)-t(\alpha)]<0\nonumber\\
&\Leftrightarrow&\int_{\alpha+1-{(\alpha+1)}^{\frac{1}{2}}}^{\alpha+1+(\alpha+1)^{\frac{1}{2}}}y^{\alpha}e^{-y}dy-\alpha\int_{0}^{\alpha+\alpha^{\frac{1}{2}}}y^{\alpha-1}e^{-y}dy<0\nonumber\\
&\Leftrightarrow&\int_{\alpha+1-{(\alpha+1)}^{\frac{1}{2}}}^{\alpha+1+(\alpha+1)^{\frac{1}{2}}}y^{\alpha}e^{-y}dy-\int_{0}^{\alpha+\alpha^{\frac{1}{2}}}y^{\alpha}e^{-y}dy<(\alpha+\alpha^{\frac{1}{2}})^{\alpha}e^{-\alpha-\alpha^{\frac{1}{2}}}\nonumber\\
&\Leftrightarrow&\int_{\alpha+\alpha^{\frac{1}{2}}}^{\alpha+1+(\alpha+1)^{\frac{1}{2}}}y^{\alpha}e^{-y}dy<\int_{0}^{\alpha+1-{(\alpha+1)}^{\frac{1}{2}}}y^{\alpha}e^{-y}dy+(\alpha+\alpha^{\frac{1}{2}})^{\alpha}e^{-\alpha-\alpha^{\frac{1}{2}}}\nonumber\\
&\Leftarrow&\int_{\alpha+\alpha^{\frac{1}{2}}}^{\alpha+1+(\alpha+1)^{\frac{1}{2}}}y^{\alpha}e^{-y}dy<(\alpha+\alpha^{\frac{1}{2}})^{\alpha}e^{-\alpha-\alpha^{\frac{1}{2}}}\nonumber\\
&\Leftrightarrow&\int_{0}^{1+(\alpha+1)^{\frac{1}{2}}-\alpha^{\frac{1}{2}}}\left(1+\frac{w}{\alpha+\alpha^{\frac{1}{2}}}\right)^{\alpha}e^{-w}dw<1\nonumber\\
&\Leftarrow&\int_{0}^{1+(\alpha+1)^{\frac{1}{2}}-\alpha^{\frac{1}{2}}}\left(1+\frac{\alpha w}{\alpha+\alpha^{\frac{1}{2}}}\right)e^{-w}dw<1\nonumber\\
&\Leftrightarrow&\frac{\alpha}{\alpha+\alpha^{\frac{1}{2}}}-\frac{3\alpha+\alpha(\alpha+1)^{\frac{1}{2}}+(1-\alpha)\alpha^{\frac{1}{2}}}{\alpha+\alpha^{\frac{1}{2}}}\cdot e^{-1-(\alpha+1)^{\frac{1}{2}}+\alpha^{\frac{1}{2}}}<0\nonumber\\
&\Leftrightarrow& e^{1+(\alpha+1)^{\frac{1}{2}}-\alpha^{\frac{1}{2}}}<3+(\alpha+1)^{\frac{1}{2}}+(1-\alpha)\alpha^{-\frac{1}{2}}\nonumber\\
&\Leftrightarrow& e^{1+(\alpha+1)^{\frac{1}{2}}-\alpha^{\frac{1}{2}}}-[1+(\alpha+1)^{\frac{1}{2}}-\alpha^{\frac{1}{2}}]<2+\alpha^{-\frac{1}{2}}.
\end{eqnarray}
Set
$$
w:=(\alpha+1)^{\frac{1}{2}}-\alpha^{\frac{1}{2}}.
$$
Then,
$$
0<w<1\quad \mbox{and} \quad \alpha^{-\frac{1}{2}}=\frac{2w}{1-w^2}.
$$
Hence
\begin{eqnarray}\label{Mar26n}
&&2+\alpha^{-\frac{1}{2}}-\left\{e^{1+(\alpha+1)^{\frac{1}{2}}-\alpha^{\frac{1}{2}}}-
\left[1+(\alpha+1)^{\frac{1}{2}}-\alpha^{\frac{1}{2}}\right]\right\}\nonumber\\
&=&\frac{3+3w-3w^2-w^3}{1-w^2}-e^{1+w}\nonumber\\
&>&\frac{3+3w-3w^2-w^3}{1-w^2}-\left[\sum_{n=0}^4\frac{(1+w)^n}{n!}+\frac{2^5}{5!}\right]\nonumber\\
&=&\frac{3 + 40 w - 153 w^2 + 160 w^3 + 145 w^4 + 40 w^5 + 5 w^6}{120(1-w^2)} .
\end{eqnarray}

Define
$$
I:=3 + 40 w - 153 w^2 + 160 w^3 + 145 w^4 + 40 w^5 + 5 w^6.
$$
Set
\begin{eqnarray*}
w:=\frac{1}{1+q^2}.
\end{eqnarray*}
Then, we have
\begin{eqnarray*}
(1 + q^2)^6I=240 + 416 q^2 + 152 q^4 + 8 q^6 + 92 q^8 + 58 q^{10} + 3 q^{12}>0,
\end{eqnarray*}
which together with (\ref{Mar26m}) and  (\ref{Mar26n}) implies that (\ref{Mar18a}) holds for the case $0<\alpha\le1$.

Next we consider the case $\alpha>1$. We have that
\begin{eqnarray*}
&&t(\alpha+1)<t(\alpha)\nonumber\\
&\Leftrightarrow&\Gamma(\alpha+1)[t(\alpha+1)-t(\alpha)]<0\nonumber\\
&\Leftrightarrow&\int_{\alpha+1-{(\alpha+1)}^{\frac{1}{2}}}^{\alpha+1+(\alpha+1)^{\frac{1}{2}}}y^{\alpha}e^{-y}dy-\alpha\int_{\alpha-\alpha^{\frac{1}{2}}}^{\alpha+\alpha^{\frac{1}{2}}}y^{\alpha-1}e^{-y}dy<0\nonumber\\
&\Leftrightarrow&\int_{\alpha+1-{(\alpha+1)}^{\frac{1}{2}}}^{\alpha+1+(\alpha+1)^{\frac{1}{2}}}y^{\alpha}e^{-y}dy-\int_{\alpha-\alpha^{\frac{1}{2}}}^{\alpha+\alpha^{\frac{1}{2}}}y^{\alpha}e^{-y}dy<(\alpha+\alpha^{\frac{1}{2}})^{\alpha}e^{-\alpha-\alpha^{\frac{1}{2}}}-(\alpha-\alpha^{\frac{1}{2}})^{\alpha}e^{-\alpha+\alpha^{\frac{1}{2}}}\nonumber\\
&\Leftrightarrow&\ \ \int_{\alpha+\alpha^{\frac{1}{2}}}^{\alpha+1+(\alpha+1)^{\frac{1}{2}}}y^{\alpha}e^{-y}dy-(\alpha+\alpha^{\frac{1}{2}})^{\alpha}e^{-\alpha-\alpha^{\frac{1}{2}}}\\
&&<\int_{\alpha-\alpha^{\frac{1}{2}}}^{\alpha+1-{(\alpha+1)}^{\frac{1}{2}}}y^{\alpha}e^{-y}dy-(\alpha-\alpha^{\frac{1}{2}})^{\alpha}e^{-\alpha+\alpha^{\frac{1}{2}}}\nonumber\\
&\Leftrightarrow&\ \ (\alpha+\alpha^{\frac{1}{2}})^{\alpha}e^{-\alpha-\alpha^{\frac{1}{2}}}\left[\int_{0}^{1+(\alpha+1)^{\frac{1}{2}}-\alpha^{\frac{1}{2}}}\left(1+\frac{y}{\alpha+\alpha^{\frac{1}{2}}}\right)^{\alpha}e^{-y}dy-1\right]\nonumber\\
&&<(\alpha-\alpha^{\frac{1}{2}})^{\alpha}e^{-\alpha+\alpha^{\frac{1}{2}}}\left[\int_{0}^{1-(\alpha+1)^{\frac{1}{2}}+\alpha^{\frac{1}{2}}}\left(1+\frac{y}{\alpha-\alpha^{\frac{1}{2}}}\right)^{\alpha}e^{-y}dy-1\right].
\end{eqnarray*}
The remainder of this section is devoted to proving the following inequality:
\begin{eqnarray}\label{Mar17}
\int_{0}^{1+(\alpha+1)^{\frac{1}{2}}-\alpha^{\frac{1}{2}}}\left(1+\frac{y}{\alpha+\alpha^{\frac{1}{2}}}\right)^{\alpha}e^{-y}dy<1<\int_{0}^{1-(\alpha+1)^{\frac{1}{2}}+\alpha^{\frac{1}{2}}}\left(1+\frac{y}{\alpha-\alpha^{\frac{1}{2}}}\right)^{\alpha}e^{-y}dy,\ \ \forall\alpha>1,\ \
\end{eqnarray}
which implies (\ref{Mar18a}). It is a bit surprising that inequality (\ref{Mar17}) is very delicate, which seems to be unknown in the literature.

\subsection{Proof of the ``$<1$" part of  inequality (\ref{Mar17})}

Note that
\begin{eqnarray*}
&&\frac{d}{dy}\left[\left(1+\frac{y}{\alpha+\alpha^{\frac{1}{2}}}\right)^{\alpha}e^{-y}\right]'\\
&=&\frac{\alpha}{\alpha+\alpha^{\frac{1}{2}}}\left(1+\frac{y}{\alpha+\alpha^{\frac{1}{2}}}\right)^{\alpha-1}e^{-y}-\left(1+\frac{y}{\alpha+\alpha^{\frac{1}{2}}}\right)^{\alpha}e^{-y}\\
&=&\frac{-\left(\alpha^{\frac{1}{2}}+y\right)}
{\alpha+\alpha^{\frac{1}{2}}}
\left(1+\frac{y}{\alpha+\alpha^{\frac{1}{2}}}\right)^{\alpha-1}e^{-y}\\
&<&0.
\end{eqnarray*}
Then, we have
\begin{eqnarray*}
&&\int_{0}^{1+(\alpha+1)^{\frac{1}{2}}-\alpha^{\frac{1}{2}}}\left(1+\frac{y}{\alpha+\alpha^{\frac{1}{2}}}\right)^{\alpha}e^{-y}dy\\
&<&\int_{0}^{1}\left(1+\frac{y}{\alpha+\alpha^{\frac{1}{2}}}\right)^{\alpha}e^{-y}dy+
\left[(\alpha+1)^{\frac{1}{2}}-\alpha^{\frac{1}{2}}\right]\left(1+\frac{1}{\alpha+\alpha^{\frac{1}{2}}}\right)^{\alpha}e^{-1}.
\end{eqnarray*}
Hence, to prove the ``$<1$" part of  inequality (\ref{Mar17}), it suffices to show that
$$
\left[(\alpha+1)^{\frac{1}{2}}-\alpha^{\frac{1}{2}}\right]\left(1+\frac{1}{\alpha+\alpha^{\frac{1}{2}}}\right)^{\alpha}e^{-1}<\int_{0}^{1}\left[1-\left(1+\frac{y}{\alpha+\alpha^{\frac{1}{2}}}\right)^{\alpha}e^{-y}\right]dy.
$$

We have
\begin{eqnarray*}
&&\frac{d}{dy}\left[1-\left(1+\frac{y}{\alpha+\alpha^{\frac{1}{2}}}\right)^{\alpha}e^{-y}\right]\\
&=&
-\frac{\alpha}{\alpha+\alpha^{\frac{1}{2}}}
\left(1+\frac{y}{\alpha+\alpha^{\frac{1}{2}}}\right)^{\alpha-1}e^{-y}
+\left(1+\frac{y}{\alpha+\alpha^{\frac{1}{2}}}\right)^{\alpha}e^{-y}\\
&=&\frac{\alpha^{\frac{1}{2}}+y}{\alpha+\alpha^{\frac{1}{2}}}
\left(1+\frac{y}{\alpha+\alpha^{\frac{1}{2}}}\right)^{\alpha-1}e^{-y}\\
&>&0,
\end{eqnarray*}
and
\begin{eqnarray*}
&&\frac{d^2}{dy^2}\left[1-\left(1+\frac{y}{\alpha+\alpha^{\frac{1}{2}}}\right)^{\alpha}e^{-y}\right]\\
&=&\frac{2\alpha}{\alpha+\alpha^{\frac{1}{2}}}\left(1+\frac{y}{\alpha+\alpha^{\frac{1}{2}}}\right)^{\alpha-1}e^{-y}
-\frac{\alpha(\alpha-1)}{(\alpha+\alpha^{\frac{1}{2}})^2}\left(1+\frac{y}{\alpha+\alpha^{\frac{1}{2}}}\right)^{\alpha-2}e^{-y}-\left(1+\frac{y}{\alpha+\alpha^{\frac{1}{2}}}\right)^{\alpha}e^{-y}\\
&=&\frac{-y(2\alpha^{\frac{1}{2}}+y)}{(\alpha+\alpha^{\frac{1}{2}})^2}\left(1+\frac{y}{\alpha+\alpha^{\frac{1}{2}}}\right)^{\alpha-2}e^{-y}\\
&<&0.
\end{eqnarray*}
Then, for fixed $\alpha$,
$1-\left(1+\frac{y}{\alpha+\alpha^{\frac{1}{2}}}\right)^{\alpha}e^{-y}$ is an increasing and concave function of $y$ on $[0,1]$. Hence, we have that
\begin{eqnarray*}
&&\int_{0}^{1}\left[1-\left(1+\frac{y}{\alpha+\alpha^{\frac{1}{2}}}\right)^{\alpha}e^{-y}\right]dy\\
&>&\frac{2}{4}\left[1-\left(1+\frac{1/2}{\alpha+\alpha^{\frac{1}{2}}}\right)^{\alpha}e^{-1/2}\right]+\frac{1}{4}\left[1-\left(1+\frac{1}{\alpha+\alpha^{\frac{1}{2}}}\right)^{\alpha}e^{-1}\right].
\end{eqnarray*}
Thus,  to complete the proof of the ``$<1$" part of  inequality (\ref{Mar17}), it suffices to show that
\begin{eqnarray*}
&&\left[(\alpha+1)^{\frac{1}{2}}-\alpha^{\frac{1}{2}}\right]\left(1+\frac{1}{\alpha+\alpha^{\frac{1}{2}}}\right)^{\alpha}e^{-1}\\
&<&\frac{2}{4}\left[1-\left(1+\frac{1/2}{\alpha+\alpha^{\frac{1}{2}}}\right)^{\alpha}e^{-1/2}\right]+\frac{1}{4}\left[1-\left(1+\frac{1}{\alpha+\alpha^{\frac{1}{2}}}\right)^{\alpha}e^{-1}\right],
\end{eqnarray*}
which is equivalent to
\begin{eqnarray}\label{R1}
\left\{1+4\left[(\alpha+1)^{\frac{1}{2}}-\alpha^{\frac{1}{2}}\right]\right\}\left(1+\frac{1}{\alpha+\alpha^{\frac{1}{2}}}\right)^{\alpha}e^{-1}+2\left(1+\frac{1/2}{\alpha+\alpha^{\frac{1}{2}}}\right)^{\alpha}e^{-1/2}<3,\ \ \ \ \forall\alpha>1.
\end{eqnarray}
The proof of (\ref{R1}) will be given in Section 4.

\subsection{Proof of the ``$>1$" part of  inequality (\ref{Mar17})}

Note that
$$
1-4[(\alpha+1)^{\frac{1}{2}}-\alpha^{\frac{1}{2}}]=0\Leftrightarrow \alpha=\left(\frac{15}{8}\right)^2.
$$
We consider the three cases $1<\alpha\le 2$, $2<\alpha<\left(\frac{15}{8}\right)^2$ and $\alpha\ge \left(\frac{15}{8}\right)^2$ separately.

\noindent {\bf  Case 1:}  $1<\alpha\le 2$.\ \
By Taylor's formula, we have
\begin{eqnarray*}
&&\int_{0}^{1-(\alpha+1)^{\frac{1}{2}}+\alpha^{\frac{1}{2}}}\left(1+\frac{y}{\alpha-\alpha^{\frac{1}{2}}}\right)^{\alpha}e^{-y}dy-1\\
&>&\int_{0}^{1-(\alpha+1)^{\frac{1}{2}}+\alpha^{\frac{1}{2}}}\left(1+\frac{\alpha y}{\alpha-\alpha^{\frac{1}{2}}}\right)e^{-y}dy-1\\
&=&\frac{\alpha}{\alpha-\alpha^{\frac{1}{2}}}-\frac{3\alpha+(\alpha-1)\alpha^{\frac{1}{2}}
-\alpha(\alpha+1)^{\frac{1}{2}}}{\alpha-\alpha^{\frac{1}{2}}}\cdot e^{-\left[1-(\alpha+1)^{\frac{1}{2}}+\alpha^{\frac{1}{2}}\right]}\\
&=&\frac{\alpha e^{-\left[1-(\alpha+1)^{\frac{1}{2}}+\alpha^{\frac{1}{2}}\right]}}{\alpha-\alpha^{\frac{1}{2}}}\left[e^{1-(\alpha+1)^{\frac{1}{2}}+\alpha^{\frac{1}{2}}}-3-(\alpha-1)\alpha^{-\frac{1}{2}}+(\alpha+1)^{\frac{1}{2}}\right]\\
&=&\frac{\alpha e^{-\left[1-(\alpha+1)^{\frac{1}{2}}+\alpha^{\frac{1}{2}}\right]}}{\alpha-\alpha^{\frac{1}{2}}}\left\{e^{1-(\alpha+1)^{\frac{1}{2}}+\alpha^{\frac{1}{2}}}
+[(\alpha+1)^{\frac{1}{2}}-\alpha^{\frac{1}{2}}]-3+\alpha^{-\frac{1}{2}}\right\}.
\end{eqnarray*}
Set
$$
w:=(\alpha+1)^{\frac{1}{2}}-\alpha^{\frac{1}{2}}.
$$
Then,
$$
\xi:=\frac{1}{2^{\frac{1}{2}}+3^{\frac{1}{2}}}\le w<\frac{1}{1+2^{\frac{1}{2}}},
$$
and
$$\alpha^{-\frac{1}{2}}=\frac{2w}{1-w^2}.
$$


We have
\begin{eqnarray*}
e^{1-(\alpha+1)^{\frac{1}{2}}+\alpha^{\frac{1}{2}}}
+[(\alpha+1)^{\frac{1}{2}}-\alpha^{\frac{1}{2}}]-3+\alpha^{-\frac{1}{2}}
=e^{1-w}+w-3+\frac{2w}{1-w^2},
\end{eqnarray*}
and
\begin{eqnarray*}
\left(e^{1-w}+w-3+\frac{2w}{1-w^2}\right)'
&=&-e^{1-w}+1+\frac{2(1+w^2)}{(1-w^2)^2}\\
&\ge&-e^{1-\xi}+1+\frac{2(1+\xi^2)}{(1-\xi^2)^2}\\
&=&1.746594.
\end{eqnarray*}
It follows that
\begin{eqnarray*}
&&e^{1-(\alpha+1)^{\frac{1}{2}}+\alpha^{\frac{1}{2}}}
+[(\alpha+1)^{\frac{1}{2}}-\alpha^{\frac{1}{2}}]-3+\alpha^{-\frac{1}{2}}\nonumber\\
&\ge& e^{1-\xi}+\xi-3+\frac{2\xi}{1-\xi^2}\\
&=&0.003095392.
\end{eqnarray*}
Thus, the ``$>1$" part of  inequality (\ref{Mar17}) holds for $1<\alpha\le 2$.

\noindent {\bf  Case 2:}  $2<\alpha< \left(\frac{15}{8}\right)^2$.\ \
By Taylor's formula, we have
\begin{eqnarray*}
&&\int_{0}^{1-(\alpha+1)^{\frac{1}{2}}+\alpha^{\frac{1}{2}}}\left(1+\frac{y}{\alpha-\alpha^{\frac{1}{2}}}\right)^{\alpha}e^{-y}dy-1\\
&>&\int_{0}^{1-(\alpha+1)^{\frac{1}{2}}+\alpha^{\frac{1}{2}}}\left[1+\frac{\alpha y}{\alpha-\alpha^{\frac{1}{2}}}+\frac{\alpha(\alpha-1) y^2}{2(\alpha-\alpha^{\frac{1}{2}})^2}\right]e^{-y}dy-1\\
&=&\frac{2\alpha^{\frac{1}{2}}+1}{\alpha^{\frac{1}{2}}-1}-\frac{\alpha\alpha^{\frac{1}{2}}-\alpha(\alpha+1)^{\frac{1}{2}}+4\alpha-4[\alpha(\alpha+1)]^{\frac{1}{2}}+8\alpha^{\frac{1}{2}}-2(\alpha+1)^{\frac{1}{2}}+2}{\alpha^{\frac{1}{2}}-1}\cdot e^{-\left[1-(\alpha+1)^{\frac{1}{2}}+\alpha^{\frac{1}{2}}\right]}\\
&=&\frac{e^{-\left[1-(\alpha+1)^{\frac{1}{2}}+\alpha^{\frac{1}{2}}\right]}}
{\alpha^{\frac{1}{2}-1}}\\
&&\cdot\left\{(2\alpha^{\frac{1}{2}}+1)e^{1-(\alpha+1)^{\frac{1}{2}}+\alpha^{\frac{1}{2}}}-
\left[\alpha\alpha^{\frac{1}{2}}-\alpha(\alpha+1)^{\frac{1}{2}}+4\alpha-
4[\alpha(\alpha+1)]^{\frac{1}{2}}+8\alpha^{\frac{1}{2}}-2(\alpha+1)^{\frac{1}{2}}+2\right]\right\}\\
&>&\frac{e^{-[1-(\alpha+1)^{\frac{1}{2}}+\alpha^{\frac{1}{2}}]}}
{\alpha^{\frac{1}{2}-1}}\Bigg\{(2\alpha^{\frac{1}{2}}+1)
\sum_{n=0}^4\frac{\left[1-(\alpha+1)^{\frac{1}{2}}+\alpha^{\frac{1}{2}}\right]^n}{n!}\\
&&\ \ \ \ \ \ \ \ -\left[\alpha\alpha^{\frac{1}{2}}-\alpha(\alpha+1)^{\frac{1}{2}}+4\alpha-
4[\alpha(\alpha+1)]^{\frac{1}{2}}+8\alpha^{\frac{1}{2}}-2(\alpha+1)^{\frac{1}{2}}+2\right]\Bigg\}.
\end{eqnarray*}

Define
\begin{eqnarray*}
J&:=&(2\alpha^{\frac{1}{2}}+1)\sum_{n=0}^4
\frac{\left[1-(\alpha+1)^{\frac{1}{2}}+\alpha^{\frac{1}{2}}\right]^n}{n!}\\
&&-\left[\alpha\alpha^{\frac{1}{2}}-\alpha(\alpha+1)^{\frac{1}{2}}+4\alpha-
4[\alpha(\alpha+1)]^{\frac{1}{2}}+8\alpha^{\frac{1}{2}}-2(\alpha+1)^{\frac{1}{2}}+2\right].
\end{eqnarray*}
Set
$$
w:=(\alpha+1)^{\frac{1}{2}}-\alpha^{\frac{1}{2}}.
$$
Then,
$$
\frac{1}{4}< w<\frac{1}{2^{\frac{1}{2}}+3^{\frac{1}{2}}},\ \ \ \ \alpha^{\frac{1}{2}}=\frac{1-w^2}{2w},\ \ \ \ (\alpha+1)^{\frac{1}{2}}=\frac{1+w^2}{2w}.
$$
Hence
\begin{eqnarray*}
J&=&\frac{1+w-w^2}{w}\sum_{n=0}^4\frac{(1-w)^n}{n!}+\frac{w^4-8w^3+18w^2-11}{4w}\\
&=&\frac{-1 + w + 9 w^2 + 38 w^3 - 31 w^4 + 9 w^5 - w^6}{24 w}\\
&>&\frac{-1 + \frac{1}{4} + (9) (\frac{1}{4} )^2 + (7)(\frac{1}{4})^3+ (31) (\frac{1}{4})^3(\frac{1}{2}) +(8)(\frac{1}{4})^5 +(\frac{1}{4})^5(\frac{1}{2})}{24 w}\\
&=&\frac{0.1723633}{24w}\\
&>&0.
\end{eqnarray*}
Thus, the ``$>1$" part of   inequality (\ref{Mar17}) holds for $2<\alpha< \left(\frac{15}{8}\right)^2$.

\noindent {\bf  Case 3:}  $\alpha\ge \left(\frac{15}{8}\right)^2$.\ \
Note that for $0\leq y\leq 1-(\alpha+1)^{\frac{1}{2}}+\alpha^{\frac{1}{2}}$, we have
\begin{eqnarray*}
&&\frac{d}{dy}\left[\left(1+\frac{y}{\alpha-\alpha^{\frac{1}{2}}}\right)^{\alpha}e^{-y}\right]'\\
&=&\frac{\alpha}{\alpha-\alpha^{\frac{1}{2}}}\left(1+\frac{y}{\alpha-\alpha^{\frac{1}{2}}}\right)^{\alpha-1}e^{-y}-\left(1+\frac{y}{\alpha-\alpha^{\frac{1}{2}}}\right)^{\alpha}e^{-y}\\
&=&\frac{\alpha^{\frac{1}{2}}-y}{\alpha-\alpha^{\frac{1}{2}}}\left(1+\frac{y}{\alpha-\alpha^{\frac{1}{2}}}\right)^{\alpha-1}e^{-y}\\
&>&0.
\end{eqnarray*}
Then,
\begin{eqnarray*}
&&\int_{0}^{1-(\alpha+1)^{\frac{1}{2}}+\alpha^{\frac{1}{2}}}\left(1+\frac{y}{\alpha-\alpha^{\frac{1}{2}}}\right)^{\alpha}e^{-y}dy\\
&>&\int_{0}^{1}\left(1+\frac{y}{\alpha-\alpha^{\frac{1}{2}}}\right)^{\alpha}e^{-y}dy-
\left[(\alpha+1)^{\frac{1}{2}}-\alpha^{\frac{1}{2}}\right]\left(1+\frac{1}{\alpha-\alpha^{\frac{1}{2}}}\right)^{\alpha}e^{-1}.
\end{eqnarray*}
Hence, to prove the ``$>1$" part of  inequality (\ref{Mar17}), it suffices to show that
$$
\left[(\alpha+1)^{\frac{1}{2}}-\alpha^{\frac{1}{2}}\right]
\left(1+\frac{1}{\alpha-\alpha^{\frac{1}{2}}}\right)^{\alpha}e^{-1}<
\int_{0}^{1}\left[\left(1+\frac{y}{\alpha-\alpha^{\frac{1}{2}}}\right)^{\alpha}e^{-y}-1\right]dy.
$$

For $y\in [0,1]$, we have
\begin{eqnarray*}
&&\frac{d}{dy}\left[\left(1+\frac{y}{\alpha-\alpha^{\frac{1}{2}}}\right)^{\alpha}e^{-y}-1\right]\\
&=&
\frac{\alpha}{\alpha-\alpha^{\frac{1}{2}}}
\left(1+\frac{y}{\alpha-\alpha^{\frac{1}{2}}}\right)^{\alpha-1}e^{-y}-
\left(1+\frac{y}{\alpha-\alpha^{\frac{1}{2}}}\right)^{\alpha}e^{-y}\\
&=&\frac{\alpha^{\frac{1}{2}}-y}{\alpha-\alpha^{\frac{1}{2}}} \left(1+\frac{y}{\alpha-\alpha^{\frac{1}{2}}}\right)^{\alpha-1}e^{-y}\\
&>&0,
\end{eqnarray*}
and
\begin{eqnarray*}
&&\frac{d^2}{dy^2}\left[\left(1+\frac{y}{\alpha-\alpha^{\frac{1}{2}}}\right)^{\alpha}e^{-y}-1\right]\\
&=&-\frac{2\alpha}{\alpha-\alpha^{\frac{1}{2}}}\left(1+\frac{y}{\alpha-\alpha^{\frac{1}{2}}}\right)^{\alpha-1}e^{-y}
+\frac{\alpha(\alpha-1)}{(\alpha-\alpha^{\frac{1}{2}})^2}\left(1+\frac{y}{\alpha-\alpha^{\frac{1}{2}}}\right)^{\alpha-2}e^{-y}+\left(1+\frac{y}{\alpha-\alpha^{\frac{1}{2}}}\right)^{\alpha}e^{-y}\\
&=&\frac{-y(2\alpha^{\frac{1}{2}}-y)}{(\alpha-\alpha^{\frac{1}{2}})^2}\left(1+\frac{y}{\alpha-\alpha^{\frac{1}{2}}}\right)^{\alpha-2}e^{-y}\\
&<&0.
\end{eqnarray*}
Then,  for fixed $\alpha$,
$\left(1+\frac{y}{\alpha-\alpha^{\frac{1}{2}}}\right)^{\alpha}e^{-y}-1$ is an increasing and concave function  of $y$ on $[0,1]$. Hence, we have that
\begin{eqnarray*}
&&\int_{0}^{1}\left[\left(1+\frac{y}{\alpha-\alpha^{\frac{1}{2}}}\right)^{\alpha}e^{-y}-1\right]dy\\
&>&\frac{2}{4}\left[\left(1+\frac{1/2}{\alpha-\alpha^{\frac{1}{2}}}\right)^{\alpha}e^{-1/2}-1\right]+\frac{1}{4}\left[\left(1+\frac{1}{\alpha-\alpha^{\frac{1}{2}}}\right)^{\alpha}e^{-1}-1\right].
\end{eqnarray*}
Thus,  to complete the proof of the ``$>1$" part of  inequality (\ref{Mar17}), it suffices to show that
\begin{eqnarray*}
&&\left[(\alpha+1)^{\frac{1}{2}}-\alpha^{\frac{1}{2}}\right]
\left(1+\frac{1}{\alpha-\alpha^{\frac{1}{2}}}\right)^{\alpha}e^{-1}\\
&<&\frac{2}{4}\left[\left(1+\frac{1/2}{\alpha-\alpha^{\frac{1}{2}}}\right)^{\alpha}e^{-1/2}-1\right]+\frac{1}{4}\left[\left(1+\frac{1}{\alpha-\alpha^{\frac{1}{2}}}\right)^{\alpha}e^{-1}-1\right],
\end{eqnarray*}
which is equivalent to
\begin{eqnarray}\label{R2}
\left\{1-4[(\alpha+1)^{\frac{1}{2}}-\alpha^{\frac{1}{2}}]\right\}
\left(1+\frac{1}{\alpha-\alpha^{\frac{1}{2}}}\right)^{\alpha}e^{-1}+2
\left(1+\frac{1/2}{\alpha-\alpha^{\frac{1}{2}}}\right)^{\alpha}e^{-1/2}>3,\ \ \ \forall \alpha\ge\left(\frac{15}{8}\right)^2.\ \ \ \
\end{eqnarray}
The proof of (\ref{R2}) will be given in Section 4.

\section{Proofs of  inequalities (\ref{R1}) and (\ref{R2})}\setcounter{equation}{0}

\subsection{Proof of  (\ref{R1})}

Define
$$
\tau_+:=\frac{1}{\alpha+\alpha^{\frac{1}{2}}},\ \ \ \ \eta_+:=\frac{\tau_+}{2}.
$$
We have
\begin{eqnarray}\label{Mar25a}
&&\left\{1+4[(\alpha+1)^{\frac{1}{2}}-\alpha^{\frac{1}{2}}]\right\}
\left(1+\frac{1}{\alpha+\alpha^{\frac{1}{2}}}\right)^{\alpha}e^{-1}+
2\left(1+\frac{1/2}{\alpha+\alpha^{\frac{1}{2}}}\right)^{\alpha}e^{-1/2}\nonumber\\
&=&\left\{1+4[(\alpha+1)^{\frac{1}{2}}-\alpha^{\frac{1}{2}}]\right\}
e^{-1+\alpha\ln(1+\tau_+)}+2e^{-\frac{1}{2}+\alpha\ln(1+\eta_+)}\nonumber\\
&<&\left\{1+4[(\alpha+1)^{\frac{1}{2}}-\alpha^{\frac{1}{2}}]\right\}
e^{-1+\alpha\left(\tau_+-\frac{\tau_+^2}{2}+\frac{\tau_+^3}{3}-\frac{\tau_+^4}{4}
+\frac{\tau_+^5}{5}\right)}+2e^{-\frac{1}{2}+\alpha\left(\eta_+-\frac{\eta_+^2}{2}+
\frac{\eta_+^3}{3}-\frac{\eta_+^4}{4}+\frac{\eta_+^5}{5}\right)}.\ \ \ \
\end{eqnarray}
Set
$$
w:=(\alpha+1)^{\frac{1}{2}}-\alpha^{\frac{1}{2}}.
$$
Then, by the condition $\alpha>1$, we get
$$
0< w<\frac{1}{2},\ \ \ \ 1-w^2>0,\ \ \ \ 1+2w-w^2>0,
$$
and
$$
\alpha=\frac{(1-w^2)^2}{4w^2},\ \ \ \ \tau_+=\frac{4w^2}{(1-w^2)(1+2w-w^2)},\ \ \ \ \eta_+=\frac{2w^2}{(1-w^2)(1+2w-w^2)}.
$$

Define
$$
P_+:=-1+\alpha\left(\tau_+-\frac{\tau_+^2}{2}+\frac{\tau_+^3}{3}-\frac{\tau_+^4}{4}+\frac{\tau_+^5}{5}\right),\ \ \ \ Q_+:=-\frac{1}{2}+\alpha\left(\eta_+-\frac{\eta_+^2}{2}+\frac{\eta_+^3}{3}-\frac{\eta_+^4}{4}+\frac{\eta_+^5}{5}\right).
$$
By virtue of {\bf Mathematica}, we get
\begin{eqnarray*}
F_+&:=&15(1 - w^2)^3 (1+2 w - w^2)^5P_+\\
&=&2 w (-15 - 135 w - 345 w^2 + 190 w^3 + 1735 w^4 + 495 w^5 -
   3615 w^6- 716 w^7 \\
&&\ \  \ \  + 3615 w^8+ 495 w^9 - 1735 w^{10} + 190 w^{11} +
   345 w^{12} - 135 w^{13} + 15 w^{14}).
\end{eqnarray*}
Set
\begin{eqnarray}\label{T1}
w:=\frac{1}{2(1+q^2)}.
\end{eqnarray}
We have
\begin{eqnarray*}
G_+&:=&16384 (1 + q^2)^{14}\cdot\frac{F_+}{2w}\\
&=&-1140603 - 17129046 q^2 - 115786348 q^4 - 468301840 q^6 -
 1267262160 q^8 - 2427446688 q^{10}\\
&& - 3393664576 q^{12} -
 3517163008 q^{14} - 2715321600 q^{16} - 1554209280 q^{18} -
 649507840 q^{20}\\
&& - 192286720 q^{22} - 38154240 q^{24} - 4546560 q^{26} -
 245760 q^{28},
\end{eqnarray*}
which implies that $P_+$ is negative.  By virtue of {\bf Mathematica}, we get
\begin{eqnarray*}
H_+&:=&30 (1 - w^2)^3 (1 + 2 w - w^2)^5Q_+\\
&=&w (-30 - 255 w - 600 w^2 + 410 w^3 + 2900 w^4 + 705 w^5 - 5550 w^6 -
   1672 w^7   \\
&&\ \  \ \ + 5550 w^8+ 705 w^9 - 2900 w^{10}+ 410 w^{11} + 600 w^{12} -
   255 w^{13} + 30 w^{14}).
\end{eqnarray*}
Further, by the transformation (\ref{T1}), we obtain that
\begin{eqnarray*}
I_+&:=&8192 (1 + q^2)^{14}\cdot\frac{H_+}{w}\\
&=&-1083048 - 16069911 q^2 - 108024568 q^4 - 435858040 q^6 -
 1178745360 q^8 - 2259543408 q^{10}\\
&& - 3165284416 q^{12} -
 3291555328 q^{14} - 2553515520 q^{16} - 1471031040 q^{18} -
 619724800 q^{20}\\
&& - 185251840 q^{22} - 37171200 q^{24} - 4485120 q^{26} -
 245760 q^{28},
\end{eqnarray*}
which implies that $Q_+$ is negative.
Thus, by (\ref{Mar25a}), we get
\begin{eqnarray}\label{Mar25b}
&&\left\{1+4[(\alpha+1)^{\frac{1}{2}}-\alpha^{\frac{1}{2}}]\right\}
\left(1+\frac{1}{\alpha+\alpha^{\frac{1}{2}}}\right)^{\alpha}e^{-1}+2\left(1+\frac{1/2}{\alpha+\alpha^{\frac{1}{2}}}\right)^{\alpha}e^{-1/2}-3\nonumber\\
&<&(1+4w)e^{P_+}+2e^{Q_+}-3\nonumber\\
&<&(1+4w)\left(1+P_++\frac{P_+^2}{2}+\frac{P_+^3}{3!}+\frac{P_+^4}{4!}\right)+2\left(1+Q_++\frac{Q_+^2}{2}+\frac{Q_+^3}{3!}+\frac{Q_+^4}{4!}\right)-3\nonumber\\
&:=&R_+.
\end{eqnarray}

Define
$$
L_+:=9720000 (1 - w^2)^{12} (1 + 2 w - w^2)^{20}w^{-3}R_+,
$$
and
$$
V_+:=-18014398509481984 (1 + q^2)^{62}L_+.
$$
By virtue of {\bf Mathematica}, we get
\begin{verbatim}
V+=23565171557938261664962395 +
 1985238765536369188253388462 q^2 +
 76017937191609745093093565184 q^4 +
 1815476155917282265018752272232 q^6 +
 30868042081839055982554050213660 q^8 +
 401897536051918258546845673711320 q^10 +
 4195397709111549929883773768957292 q^12 +
 36238699732610615067411794056699104 q^14 +
 265002286089679374723172860122766982 q^16 +
 1669237124849349342077586449716389470 q^18 +
 9179934813394932229977676436328785920 q^20 +
 44555295354320392501114611345123622400 q^22 +
 192537160208281140648975165919934835200 q^24 +
 746181252269526741637909507751082171520 q^26 +
 2609372572626683435719917787491018652160 q^28 +
 8276209631283583168755734561689661224960 q^30 +
 23913569456882144063241575623509484876800 q^32 +
 63185851825755484161960668172292699909120 q^34 +
 153171842040744452342444666253152790732800 q^36 +
 341628452529444844018632398179833131991040 q^38 +
 702775058219816773204544225960017412751360 q^40 +
 1336281286191807830241756821296507838955520 q^42 +
 2352931028757956298911312671496544634142720 q^44 +
 3842851078067257537573706091171559356497920 q^46 +
 5829597689354288278514821031866786159656960 q^48 +
 8224048629268397888867021111129844469596160 q^50 +
 10800344227796355322915422778734559920914432 q^52 +
 13214925874596962589048294078754839901241344 q^54 +
 15075437985869745487585690312752934894436352 q^56 +
 16043253396277674458218536937392302561689600 q^58 +
 15933483495160554733717977505944446676500480 q^60 +
 14772181225346687498328707734409833005187072 q^62 +
 12786642631638024827680853914736629691449344 q^64 +
 10333607587858511627607426462208954269171712 q^66 +
 7796135691442288431829216828566360534548480 q^68 +
 5489455045169343972022956242977322445045760 q^70 +
 3606053976460179205452292516224406577479680 q^72 +
 2208802483012122361676530694278736018145280 q^74 +
 1260683716848402925787749700503070572544000 q^76 +
 669904634456217504368236284579198848204800 q^78 +
 331081492020125941589702200450146757509120 q^80 +
 152001230583526972614803279225239581491200 q^82 +
 64734238129983061042604032362158017740800 q^84 +
 25531661312772564369272230408940525977600 q^86 +
 9307900274880934185285303838052660019200 q^88 +
 3129622227458365558314112917099983667200 q^90 +
 968033669852811173795309928049750835200 q^92 +
 274640462267821497605095290057418342400 q^94 +
 71224038857339584104332613373237657600 q^96 +
 16816854347488682979216179348688076800 q^98 +
 3598246400042953802386878316412928000 q^100 +
 693860503839280523254707460767744000 q^102 +
 119794916777653504670468143054848000 q^104 +
 18371891299642536871251828277248000 q^106 +
 2478647665316721166509051740160000 q^108 +
 290656583441325139861690122240000 q^110 +
 29171448597603811259616067584000 q^112 +
 2455470675466333890778497024000 q^114 +
 168582129968383522599075840000 q^116 +
 9065459012974557989437440000 q^118 +
 358068461185282678456320000 q^120 +
 9236522547766697656320000 q^122 +
 116733302341443256320000 q^124.
\end{verbatim}
Note that all terms in the above expansion are positive. Then, $R_+$ is negative. Therefore, the proof is complete by (\ref{Mar25b}).

\subsection{Proof of  (\ref{R2})}

Define
$$
\tau_-:=\frac{1}{\alpha-\alpha^{\frac{1}{2}}},\ \ \ \ \eta_-:=\frac{\tau_-}{2}.
$$
We have
\begin{eqnarray}\label{Mar25aa}
&&\left\{1-4[(\alpha+1)^{\frac{1}{2}}-\alpha^{\frac{1}{2}}]\right\}
\left(1+\frac{1}{\alpha-\alpha^{\frac{1}{2}}}\right)^{\alpha}e^{-1}+
2\left(1+\frac{1/2}{\alpha-\alpha^{\frac{1}{2}}}\right)^{\alpha}e^{-1/2}\nonumber\\
&=&\left\{1-4[(\alpha+1)^{\frac{1}{2}}-\alpha^{\frac{1}{2}}]\right\}
e^{-1+\alpha\ln(1+\tau_-)}+2e^{-\frac{1}{2}+\alpha\ln(1+\eta_-)}\nonumber\\
&>&\left\{1-4[(\alpha+1)^{\frac{1}{2}}-\alpha^{\frac{1}{2}}]\right\}
e^{-1+\alpha\left(\tau_--\frac{\tau_-^2}{2}+\frac{\tau_+^3}{3}-\frac{\tau_-^4}{4}\right)}
+2e^{-\frac{1}{2}+\alpha\left(\eta_--\frac{\eta_-^2}{2}+\frac{\eta_-^3}{3}-
\frac{\eta_-^4}{4}\right)}.
\end{eqnarray}
Set
$$
w:=(\alpha+1)^{\frac{1}{2}}-\alpha^{\frac{1}{2}}.
$$
Then, by the condition $\alpha\geq \left(\frac{15}{8}\right)^2$, we have
$$
0< w\le\frac{1}{4},\ \ \ \ 1-w^2>0,\ \ \ \ 1-2w-w^2>0,
$$
and
$$
\alpha=\frac{(1-w^2)^2}{4w^2},\ \ \ \ \tau_-=\frac{4w^2}{(1-w^2)(1-2w-w^2)},\ \ \ \ \eta_-=\frac{2w^2}{(1-w^2)(1-2w-w^2)}.
$$

Define
$$
P_-:=-1+\alpha\left(\tau_--\frac{\tau_-^2}{2}+\frac{\tau_-^3}{3}-\frac{\tau_-^4}{4}\right),\ \ \ \ Q_-:=-\frac{1}{2}+\alpha\left(\eta_--\frac{\eta_-^2}{2}+\frac{\eta_-^3}{3}-\frac{\eta_-^4}{4}\right).
$$
We have
\begin{eqnarray*}
F_-&:=&3 (1 - w^2)^2 (1 - 2 w - w^2)^4P_-\\
&=&-2 w (-3 + 21 w - 33 w^2 - 56 w^3 + 130 w^4 + 94 w^5 - 130 w^6 - 56 w^7 +33 w^8 + 21 w^9 + 3 w^{10}).
\end{eqnarray*}
Set
\begin{eqnarray}\label{T2}
w:=\frac{1}{4(1+q^2)}.
\end{eqnarray}
We get
\begin{eqnarray*}
G_-&:=&1048576 (1 + q^2)^{10}\cdot\frac{F_-}{2w}\\
&=&128409 + 2102668 q^2 + 14459888 q^4 + 56813056 q^6 + 142035456 q^8 +
 236177408 q^{10}\\
&& + 264626176 q^{12} + 197525504 q^{14} + 94175232 q^{16} +
 25952256 q^{18} + 3145728 q^{20},
\end{eqnarray*}
which implies that $P_-$ is positive. We have
\begin{eqnarray*}
H_-&:=&6 (1 - w^2)^2 (1 - 2 w - w^2)^4Q_-\\
&=&w (6 - 39 w + 54 w^2 + 100 w^3 - 200 w^4 - 128 w^5 + 200 w^6 +
   100 w^7 - 54 w^8 - 39 w^9 - 6 w^{10}).
\end{eqnarray*}
By the transformation (\ref{T2}), we get
\begin{eqnarray*}
I_-&:=&524288 (1 + q^2)^{10}\cdot\frac{H_-}{w}\\
&=&175743 + 2666962 q^2 + 17644496 q^4 + 67116160 q^6 + 162604032 q^8 +
 262406144 q^{10}\\
&&  + 286081024 q^{12} + 208437248 q^{14} + 97320960 q^{16} +
 26345472 q^{18} + 3145728 q^{20},
\end{eqnarray*}
which implies that $Q_-$ is positive.
Thus, by (\ref{Mar25aa}), we get
\begin{eqnarray}\label{Mar25bb}
&&\left\{1-4[(\alpha+1)^{\frac{1}{2}}-\alpha^{\frac{1}{2}}]\right\}
\left(1+\frac{1}{\alpha-\alpha^{\frac{1}{2}}}\right)^{\alpha}e^{-1}+
2\left(1+\frac{1/2}{\alpha-\alpha^{\frac{1}{2}}}\right)^{\alpha}e^{-1/2}-3\nonumber\\
&>&(1-4w)e^{P_-}+2e^{Q_-}-3\nonumber\\
&>&(1-4w)\left(1+P_-+\frac{P_-^2}{2}+\frac{P_-^3}{3!}\right)+2\left(1+Q_-+\frac{Q_-^2}{2}+\frac{Q_-^3}{3!}\right)-3\nonumber\\
&:=&R_-.
\end{eqnarray}

Define
$$
L_-:=-648 (1 - w^2)^6 (1 - 2 w - w^2)^{12}w^{-3}R_-,
$$
and
$$
V_-:=-9223372036854775808 (1 + q^2)^{34}L_-.
$$
By virtue of {\bf Mathematica}, we get
\begin{verbatim}
V-=1058023271132626023 + 51541890229923566472 q^2 +
 1213009372688989850064 q^4 + 18352820646596071930240 q^6 +
 200442482186879766344000 q^8 + 1682464063207304317242816 q^10 +
 11285809233594557704985856 q^12 + 62123650712872430361438720 q^14 +
 286006349074965960756670464 q^16 +
 1116988330180696358380290048 q^18 +
 3741070988530167056939876352 q^20 +
 10836728622922107883411734528 q^22 +
 27330768999389436608140804096 q^24 +
 60330629581678789398471114752 q^26 +
 117040259220868123540341129216 q^28 +
 200169684615441568277429485568 q^30 +
 302487634652646366318853881856 q^32 +
 404485176548048584076188188672 q^34 +
 478958931700928292722105647104 q^36 +
 502214198209105947113507782656 q^38 +
 465950173145939141611449483264 q^40 +
 381911302204202972478305206272 q^42 +
 275855094787796630236532047872 q^44 +
 174972859491619945161380855808 q^46 +
 97001064005371803174963249152 q^48 +
 46708058892337905269349548032 q^50 +
 19376798019028218231165812736 q^52 +
 6852048396846611541188935680 q^54 +
 2036474622748306983572471808 q^56 +
 499070059195604547617685504 q^58 + 98184545971566239935365120 q^60 +
 14905936080354118622773248 q^62 + 1639091209276985243074560 q^64 +
 116172982490204328689664 q^66 + 3984496719921263149056 q^68.
\end{verbatim}
Note that all terms in the above expansion are positive.  Then, $R_-$ is positive. Therefore, the proof is complete by (\ref{Mar25bb}).
\vskip 0.5cm
{ \noindent {\bf\large Acknowledgements}\quad  This work was supported by the National Natural Science Foundation of China (No. 12171335), the Science Development Project of Sichuan University (No. 2020SCUNL201) and the Natural Sciences and Engineering Research Council of Canada (No. 4394-2018).

\end{document}